\documentclass{annalife}

\begin{document}

\title{ABS ALGORITHMS FOR LINEAR EQUATIONS  
AND ABSPACK}
\author{\uppercase{~E. Spedicato}\istituto{
Department of Mathematics, University of Bergamo, Bergamo, Italy}
   \and \uppercase{~Z. Xia}\istituto{Department of Applied Mathematics, 
Dalian University of Technology, Dalian, China} \and \uppercase{~E. Bodon}\istituto{Department 
of Mathematics, University of Bergamo, Bergamo, Italy} 
\and \uppercase{~A. Del Popolo}\istituto{Department of 
Mathematics, University of Bergamo, Bergamo, Italy}
}
\authorhead{\lowercase{~E. Spedicato - ~Z. Xia - ~E. Bodon - 
~A. Del Popolo}}
\titlehead{\lowercase{Abs algorithms for linear equations}}
\date{}
\maketitle
\begin{abstract}
We present the main results obtained
during a  research on ABS methods in the framework of the
project Analisi Numerica e Matematica Computazionale.
\end{abstract}
\noindent
\section{The ABS algorithms}.

ABS algorithms were  introduced by Abaffy, Broyden and Spedicato (1984), 
 to solve linear equations first in the form of the
 {\it basic ABS class}, later generalized 
 as the  {\it scaled ABS class} and 
 applied  to linear least squares, nonlinear equations 
and optimization problems, see e.g. the monographs by
Abaffy and Spedicato (1989) and Zhang, Xia and Feng (1999), or the
  bibliography by Nicolai and Spedicato 
  (1997) listing over 300 ABS papers. In this paper we  
 consider some  new results obtained in the framework of
 project Analisi Numerica e Matematica Computazionale, including the 
 performance of several codes of  ABSPACK, a
 FORTRAN package  under development.

 For later reference, we recall the scaled ABS algorithm for solving the following
 determined or underdetermined linear system, where
$rank(A)$ is arbitrary and $ A^T = (a_1,\ldots,a_m) $

$$
Ax = b \ \ \ \ \  x \in {R}^n ,\ \  b \in {R}^m , \ \ m \leq n \eqno {(1)}
$$
or
$$
a_i^Tx - b_i = 0, \ \ \  i = 1,\ldots,m  \eqno {(2)}
$$

\vskip2mm
\noindent {\bf The scaled ABS algorithm}
\vskip2mm
\begin{description}
\item[(A)] Give $x_1 \in {R}^n$  arbitrary,\ $H_1 \in {R}^{n,n}$ 
nonsingular arbitrary,  $v_1 \in R^m$  arbitrary nonzero.
Set $i = 1$.
\item[(B)] Compute the residual $r_i = Ax_i-b$. If $r_i = 0$ stop ($x_i$
solves the problem.) Otherwise compute
$s_i = H_iA^Tv_i$. If $s_i \neq 0$, then  go  to (C).  If $s_i = 0$
and $\tau = v_i^Tr_i = 0$, then set $x_{i+1} = x_{i} ,\  H_{i+1} = H_{i}$
(the $i$-th equation is redundant) 
and go to (F). Otherwise stop (the system has no solution).
\item[(C)] Compute the search vector $p_i$ by
$$
p_i = H_i^Tz_i \eqno{(3)}
$$
where $z_i \in  {R}^n$ is arbitrary save for the condition
$$
v_i^TAH_i^Tz_i \neq 0 \eqno{(4)}
$$
\item[(D)] Update the estimate of the solution by
$$
x_{i+1} = x_{i} - \alpha_ip_i, \ \ \ \alpha_i = v_i^Tr_i/v_i^TAp_i.  \eqno{(5)}
$$
\item[(E)] Update the matrix $H_i$ by
$$
H_{i+1} = H_i - H_iA^Tv_iw_i^TH_i/w_i^TH_iA^Tv_i \eqno{(6)}
$$
where $w_i \in {R}^n$ is arbitrary save for the condition
$$
w_i^TH_iA^Tv_i \neq 0.   \eqno{(7)}
$$
\item[(F)] If $i = m$,  stop  ($x_{m+1}$ solves the system). Otherwise 
give  $v_{i+1} \in R^m $  arbitrary linearly
independent from $v_1, \ldots, v_i$. Increment
$i$ by one and go to ({\bf B}).
\end{description}
Matrices $H_i$, which  are generalizations of  
(oblique) projection matrices, have been named
 {\it Abaffians} at the First International Conference on ABS methods
 (Luoyang, China, 1991). There are  alternative formulations of the
 scaled ABS algorithms, e.g. using
 vectors instead of the  square matrix $H_i$, with possible
 advantages in storage and number of operations. In a next section we will 
 show how they can be used to generate  infinite iterative
 methods.

 The  choice of the parameters $H_1, \ v_i, \ ,z_i, \  w_i$ determines
 particular methods.
The {\it basic ABS class} is obtained by taking $v_i = e_i, $ as the
$i$-th unit vector in $R^m$. 

We recall some  properties of  the scaled ABS class,
assuming that $A$ has full rank.
\begin{itemize}
\item Define  $V_i = (v_1,\ldots,v_i)$ and
$W_i = (w_1,\ldots,w_i)$. Then $H_{i+1}A^TV_i = 0$ and  $H_{i+1}^TW_i = 0$.
\item The vectors $H_iA^Tv_i, \  H_i^Tw_i $ are zero 
if and only if $a_i, \ w_i$ are respectively linearly dependent from
$a_1,\ldots,a_{i-1}, \ w_1, \ldots, w_{i-1}$. 
\item Define  $P_i = (p_1,\ldots,p_i) $ and $A_i = (a_1, \ldots,a_i)$.
Then the  implicit factorization 
$V_i^TA_i^TP_i = L_i$ holds, where $L_i$ is nonsingular lower triangular.
Hence, if $m=n$, one obtains a semiexplicit 
factorization of the inverse, with $P = P_n, \ V = V_n, \ L = L_n$
$$
A^{-1} = PL^{-1}V^T.  \eqno{(8)}
$$
For several choices of  $V$ the matrix $L$ is diagonal, hence
formula (8) gives a fully explicit factorization of the inverse as a
byproduct of the ABS solution of a linear system.

\item The general solution of system (1) can be written as follows,
with $q \in R^n$ arbitrary
$$
x = x_{m+1} + H_{m+1}^Tq    \eqno{(10)}
$$
\item The Abaffian can be written in term of the parameter matrices
as
$$
H_{i+1} = H_1 - H_1A^TV_i(W_i^TH_1A^TV_i)^{-1}W_i^TH_1. \eqno{(11)}
$$
\end{itemize}
Letting $V = V_m, \ W = W_m$, one can show that the parameter matrices $H_1, \ 
V, \ W$ are admissible (i.e.  condition (7) is satisfied) iff
the matrix $Q = V^TAH_1^TW$ is strongly nonsingular (i.e. it 
is LU factorizable).
This condition can  be satisfied by  exchanges
of the columns of $V$ or $W$.
If $Q$ is strongly nonsingular and we take, as is 
done in all algorithms so far considered, $z_i = w_i$, then condition
(4) is also satisfied. Analysis of the conditions under which $Q$ is not
strongly nonsingular leads, when dealing with Krylov space 
methods in their ABS formulation, to a characterization of the topology
of the starting points that can produce a breakdown (either  a division by
zero or  a vanishing search direction) and to several ways of curing it,
 including those considered in the literature.

Two subclasses of the scaled ABS class and  particular algorithms are
now recalled.

\begin {description}
\item[(a)] The {\it conjugate direction subclass}.
This class is obtained by setting $v_i = p_i$. It 
is well defined under the condition (sufficient but not necessary)
that $A$ is symmetric and positive definite. It contains 
the ABS versions of
the  Choleski,
the Hestenes-Stiefel and the Lanczos algorithms. This class generates
all possible algorithms whose search directions are $A$-conjugate. 
If $x_1=0$, the vector
$x_{i+1}$ minimizes the energy ($A$-weighted Euclidean)
norm of the error over  $ Span(p_1,\ldots,p_i)$ and the
solution is approached monotonically from below in the energy norm.

\item[(b)] The {\it orthogonally scaled subclass}.
This class is obtained by setting $v_i = Ap_i$. It 
is well defined if $A$ has full column rank and remains 
well defined even if $m$ is greater than
$n$. It contains the ABS formulation of the QR algorithm
(the so called {\it implicit QR algorithm}),  the GMRES and 
the conjugate residual algorithms. The scaling vectors are orthogonal and
the search vectors are $A^TA$-conjugate. If $x_1=0$, 
the vector $x_{i+1}$ minimizes the
Euclidean norm of the residual over $ Span(p_1, \ldots,p_i)$ and
the solution is monotonically approached from below in the residual norm.
It can be shown that
the methods in this class can be applied to overdetermined systems 
of $ m > n$ equations, where in $n$ steps they  obtain
the solution in the least squares sense. 
\item[(c)] The {\it optimally stable subclass}. This class is obtained by setting
$v_i = A^{-T}p_i$, the inverse disappearing in the actual recursions. The
search vectors in this class are orthogonal.
If $x_1=0$, then
the vector $x_{i+1}$ is the vector of least Euclidean norm over
$ Span(p_1, \ldots,p_i)$ 
and the solution is approached
monotonically from below in the  Euclidean norm. The methods of Gram-Schmidt and
of Craig belong to this subclass. The methods in this class have minimum
error growth in the approximation to the solution according to a criterion
by Broyden.
\item[(d)] 
The {\it Huang algorithm} is obtained by the  choices $H_1 = I$, $z_i =
w_i = a_i, \ v_i = e_i$. 
A mathematically equivalent, but numerically more stable, formulation
is the so called {\it modified Huang algorithm} 
($p_i = H_i(H_ia_i)$ and
$H_{i+1} = H_i - p_ip_i^T/p_i^Tp_i$). 
Huang algorithm generates
 search vectors that are orthogonal and identical with those
obtained by the  Gram-Schmidt  procedure applied
to the rows of $A$. 
 If $x_1=0$, then  $x_{i+1}$ is the solution
with least Euclidean norm of the first $i$ equations. The solution $x^+$ with
least Euclidean norm of the whole system is approached monotonically and from below
by the sequence $x_i$. 
\item[(e)] 
The {\it implicit LU algorithm} is given by
the choices $H_1 = I,\  z_i = w_i = v_i = e_i$.
It is well defined iff $A$ is regular (i.e. all 
principal submatrices are nonsingular). Otherwise column
pivoting  has to be performed (or, if $m = n$, equation pivoting).
The Abaffian  has the following structure, with
$K_i \in R^{n-i,i}$
$$H_{i+1} = \left
	  [
	  \begin{array}{cc}
	  0      & 0 \\
	  \cdots & \cdots \\
	  0      & 0 \\
	  K_i    & I_{n-i}
	  \end{array}
	  \right ].     \eqno{(12)}
$$
implying that the matrix $P_i$ is
unit upper triangular, so that the implicit factorization $A = LP^{-1}$ is of
the LU type, with units on the diagonal.
 The algorithm requires 
 for $m  =  n$,  
\  $n^3/3$
multiplications plus lower order terms,  the same cost of
   classical LU factorization or Gaussian
elimination. 
Storage requirement  for $K_i$
 requires at most  $n^2/4$ positions, i.e.
 half the storage needed by Gaussian elimination and a fourth that
  needed by the LU factorization algorithm
(assuming that $A$ is not overwritten). Hence the implicit LU algorithm
has same arithmetic cost but uses less memory than the most efficient
classical methods.
\item[(f)] 
The {\it implicit LX algorithm}, see Spedicato, Xia and Zhang (1997), 
is defined by the choices
$H_1 = I,v_i = e_i,\ z_i = w_i = e_{k_i}$,
where $k_i$ is an integer, $1 \leq k_i \leq n$, such that
$e_{k_i}^TH_ia_i \neq 0. $ 
By a general property of the ABS
class for $A$ with full rank there is at least
one index $k_i$ such that $e_{k_i}^TH_ia_i \neq 0. $ For stability reasons
we select $k_i$ such
that $\mid e_{k_i}^TH_ia_i \mid$ is maximized. This algorithm
has the same overhead and memory requirement as the implicit LU algorithm,
but does not require  pivoting. Its computational performance is also
superior and generally better than the performance of  the classical LU
factorization algorithm with row pivoting, as available for instance in
 LAPACK or MATLAB, see Mirnia (1996). Therefore this algorithm can be
considered as {\it the most efficient 
general purpose linear solver not of the Strassen type}. The implicit LX 
algorithm
has also very important applications in a reformulation of the simplex
method for the LP problem, see Zhang, Xia and Feng (1999), where it
leads to a reduction of storage up to a factor 8 and of multiplications
up to one  order for problems where there is a small number of degrees of 
freedom, with respect to implementations based upon the classical LU
factorization.
\end{description}
\section{Solution of linear Diophantine equations}
One of our main results  has been the derivation of   
ABS methods for linear
Diophantine equations. The ABS algorithm determines 
if the Diophantine system has an
integer solution, computes a particular solution and provides a
representation of  all integer solutions. It is a generalization of
  a method proposed by Egervary (1955) for the 
 particular case of a homogeneous system.

Let ${Z}$ be the set of all integers and consider the Diophantine linear
system of equations
$$
Ax=b, \quad  \quad x\in {Z}^n,             \
A\in {Z}^{m\times n}, \ b\in {Z}^m, \ m\leq n.\eqno{(13)}
$$

It is intriguing  that while  thousands of  papers have been
written concerning nonlinear, usually polynomial, Diophantine equations
in few variables, the general linear system has attracted   much
less attention. The single  linear equation in $n$ variables was first
solved by Bertrand and Betti (1850). Egervary  was probably
the first author dealing with  
a system (albeit only the homogeneous one).  Several methods for 
the nonhomogeneous system have recently been proposed based mainly on
reduction to canonical forms.

We  recall some results from number theory.
Let $a$ and $b$ be integers. If there is integer $\gamma$ 
so that $b= \gamma a$ then
we say that $a$ divides $b$ and write $a|b$, otherwise 
we write $a \!\!\not|\ b$.
If $a_1,\dots,a_n$ are integers, not all being zero, then the greatest common
divisor ($gcd$) of these numbers is the greatest positive integer $\delta$ which
divides all
$a_i$, $i=1,\dots,n$ and we write $ \delta =gcd(a_1,\dots,a_n)$. We
note that
$\delta \geq 1$ and that 
$\delta$ can be written as an integer linear combination of
the $a_i$, i.e. $\delta = z^Ta$ for some $z\in R^n$.
One can show that $\delta$ is the least positive integer for which
the equation $a_1x_1+\dots+a_nx_n= \delta $ has an integer solution. Now
$\delta$ plays a main role in the following
\vskip3mm
	 
\noindent
{{\bf  Fundamental Theorem of the  
Linear Diophantine Equation}}
\vskip2mm
{\it
\noindent Let $a_1,\dots,a_n$ and $b$ be integer numbers. Then the
 Diophantine linear equation $a_1x_1+\dots+a_nx_n=b$
has  integer solutions if and only
if $gcd(a_1,\dots,a_n) |\ b$. 
In such a case if $n>1$ then
there is an infinite number of integer solutions. }

\vspace{0.5cm}

\noindent In order to find the general integer solution of the Diophantine
equation $a_1x_1+\dots+a_nx_n=b$, the main step is to solve
$a_1x_1+\dots+a_nx_n= \delta$,
where $\delta=gcd(a_1,\dots,a_n)$, for a special integer solution. 
There exist several
algorithms for this problem.
The basic step is the computation of $\delta$ and $z$,
often done using  the  algorithm of
Rosser (1941), which  
 avoids a too rapid growth of the intermediate integers, and which
 terminates in polynomial time, as shown
by Schrijver (1986). The scaled ABS algorithm can be applied to Diophantine  equations via a special
choice of its parameters, originating from the following considerations
and Theorems.

Suppose $x_i$ is an integer vector.
Since $x_{i+1}=x_i-\alpha_ip_i$, then $x_{i+1}$ is integer if
$\alpha_i$ and $p_i$ are integers. If $v_i^TAp_i|(v_i^Tr_i)$,
 then $\alpha_i$ is an integer. If $H_i$
and $z_i$ are respectively an integer
matrix and an integer vector, 
then $p_i=H_i^Tz_i$ is also an integer vector. Assume
$H_i$ is an integer matrix. From (6), if
$v_i^TAH_i^Tw_i$ divides all the components of $H_iA^Tv_i$, then $H_{i+1}$
is an integer matrix. 

Conditions for the existence of an integer solution and determination of all 
integer  solutions of the Diophantine system are given in the following
theorems,  generalizing the Fundamental Theorem, see 
Esmaeili, Mahdavi-Amiri and Spedicato
 (1999), or Fodor (1999)
for a different proof under somewhat less general conditions. 

\begin{Theo} \
 {\it Let $A$ be full rank and suppose that the Diophantine
system (13) is integerly solvable. Consider the  Abaffians generated
by the scaled ABS algorithm with the  parameter choices:
  $H_1$ is unimodular (i.e. both  $H_1$ and $H_1^{-1}$ are integer
  matrices);
  for $i = 1, \ldots, m$, \  $w_i$ is such that
$w_i^TH_iA^Tv_i = \delta_i,  \ \ \delta_i = gcd(H_iA^Tv_i)$.
Then the following properties are true:
\begin{description}
\item[(a)] the  Abaffians generated by the algorithm are
well-defined and are integer matrices
\item[(b)] if $y$ is a special integer solution of the first $i$
equations, then any integer solution $x$ of such equations can be
written as $x = y + H_{i+1}^Tq$ for some integer vector
$q$.
\end{description}    }
\end{Theo}
  
\begin{Theo} \
 {\it 
Let $A$ be full rank and consider the sequence of
matrices $H_i$ generated by the scaled ABS algorithm with parameter
choices as in Theorem 1. Let  $x_1$ in the scaled
ABS algorithm be an arbitrary integer vector and let $z_i$ be  such that
$z_i^TH_iA^Tv_i = gcd(H_iA^Tv_i)$.
Then system (13)         has integer solutions iff 
$gcd(H_iA^Tv_i)$ divides $v_i^Tr_i $ for $i = 1, \ldots, m$.  }
\vskip3mm
\noindent We can now state the scaled ABS algorithm for Diophantine equations.
\end{Theo}

\noindent
{\bf The  ABS Algorithm for  Diophantine Linear Equations}
\begin{enumerate}
\item[(1)] Choose $x_1\in{Z}^n$, arbitrary,  $H_1\in {Z}^{n\times n}$,
      arbitrary  unimodular. Let $i=1$.
\item[(2)] Compute $\tau_i=v_i^Tr_i$ and $s_i=H_iA^Tv_i$.
\item[(3)] {{\it If}} ($s_i=0$ and $\tau_i=0$) {{\it then}} let $x_{i+1}=x_i$,
      $H_{i+1}=H_i$, $r_{i+1}=r_i$ and {{\it go to}} step (5) (the $i$th
      equation is redundant).
      {{\it If}} ($s_i=0$ and $\tau_i\not =0$) {{\it then}} Stop (the
      $i$th equation and hence the system is incompatible).
\item[(4)] $\{s_i\not =0\}$ Compute  $\delta_i=gcd(s_i)$
      and $p_i=H_i^Tz_i$, where
      $z_i\in {Z}^n$ is an arbitrary integer vector satisfying
      $z_i^Ts_i=\delta_i$. {{\it If}} \ \
      $\delta_i\!\!\not| \tau_i$ {{\it then}} Stop 
      (the system is integerly inconsistent),
     {{\it else}}
      Compute
      $
      \alpha_i=\tau_i/\delta_i,
      $
      let
      $
      x_{i+1}=x_i-\alpha_ip_i
      $
 and update $H_i$  by
      $
      H_{i+1}=H_i-\frac{H_iA^Tv_iw_i^TH_i}{w_i^TH_iA^Tv_i}
      $
      where $w_i\in {R}^n$ is an arbitrary integer vector satisfying
       $w_i^Ts_i=\delta_i$.
\item[(5)] {{\it If}} $i=m$ {{\it then}} Stop ($x_{m+1}$ is a solution)
{{\it else}} let $i=i+1$ and {{\it go to}} step (2).
\end{enumerate}

\vspace{0.5cm}

\noindent It follows from Theorem 1 that 
if there exists a solution for the system (13), then
$x=x_{m+1}+H_{m+1}^Tq$, with arbitrary $q\in{Z}^n$, provides all
solutions of (13).  

Egervary's algorithm  for 
homogeneous Diophantine systems corresponds to the choices 
$H_1=I,\ x_1=0$ and
$w_i=z_i$, for all $i$. Egervary claimed, without  proof, that any set of
$n-m$ linearly independent rows of $H_{m+1}$ form an integer basis for the
general solution of the system. We have shown by a counterexample that
Egervary's claim is not true in general; we have also provided an
analysis of conditions under 
which $m$ rows in $H_{m+1}$ can be eliminated.

\vskip9mm

\section{The generalized implicit LU subclass} 

The  generalized implicit  LU \ (GILU) \ subclass 
is defined by taking  $v_i = z_i = w_i = e_i$ and 
 $H_1$  arbitrary nonsingular. The 
 GILU subclass is well defined iff the matrix 
$AH_1^TE_m$ is strongly nonsingular, where $E_m = (e_1,...,e_m)$.

The well definiteness condition involves the matrix 
$AH_1^T = (H_1a_1,\ldots,H_1a_m)^T$.\ $H_1^T$ can be interpreted as a right
scaling or right conditioning operator on $A$, acting in the
same way on the different rows of $A$. If $A$ is full rank but not
regular, the well definiteness condition can be satisfied
by simply taking $H_1$ as a suitable permutation matrix. By this
choice LU factorization with column pivoting is imbedded in
the GILU subclass.    It can also be shown that all
sequences $x_i$ generated by the basic ABS class can be obtained a suitable
choice of $H_1$ in the GILU subclass.

The given parameter choices
imply the following structure for the Abaffian,
$$H_{i+1} = \left
	    [
	    \begin{array}{c}
	    0_{i,n}\\
	    S_{n-i,n}
	    \end{array}
	    \right ] \eqno{(14)}
$$
where  $S_{n-i,n} \in {R}^{n-i,n}$ 
is full rank. 
The total number of multiplications is no more than 
$n^3 + O(n^2)$  for $m=n$, a substantial saving over
algorithms in the basic ABS class which may require $3n^3$ multiplications, after
the parameters are given.

The GILU subclass is related to a representation of the ABS class 
given in Abaffy and Spedicato (1989) in terms of 
$n - i$ vectors as a generalization of a method proposed by Sloboda (1978).
In such a representation one takes $w_i = z_i$ and $H_1$ arbitrary,  assuming 
that the feasible parameters $z_i$'s are known initially. Then the search 
vector can be written in the form
$p_i = u_i^i $
where $u_j^1 = H_1^{\ T}z_j,\ \ j = 1,\ldots,n$ and the vectors $u_j^i$ are updated, 
for $i = 1,\ldots,m$, by
$u_j^{i+1} = u_j^i - (a_i^Tu_j^i/a_i^Tu_i^i)u_i^i.$

The relation between the GILU subclass and the
representation in terms of $n-i$ vectors is given by the following
\begin{Theo} \
{\it Define the matrix $U_i = (u_1,\ldots,u_n)$ by $u_k = 0$ for $k < i, 
u_k = u_k^i$ for  $k \geq i$. Then 
$U_i = H_i^T ,\ \   u_i^i = p_i $
where $H_i, \ p_i$ are respectively the i-th Abaffian and search vector 
of the GILU subclass with initial Abaffian 
$H_1' = Z^TH_1, \ \ Z = (z_1,\ldots,z_n)$.}
\end{Theo}
\vskip3mm
Further analysis shows that
the matrix $H_1$ needs not be given  explicitly at the initial 
step of the algorithm. The $i$-th row of $H_1$ can be  given just at the 
$i$-th step. It can be considered as a vector parameter (right scaling 
parameter), arbitrary save that the matrix $A_i^T(H_1^T)_i$ must be 
nonsingular, where $(H_1^T)_i$ is the matrix comprising the first $i$
columns of $H_1^T$.  Therefore this formulation shows that all right preconditioners
can be imbedded in the ABS class, right preconditioning being just
equivalent to a change in the initial matrix $H_1$. 
\vskip9mm
\section{ABS Methods for KT Equations.}
\vskip5mm
The KT (Kuhn-Tucker) equations
 are the following special linear system, 
 related to the optimality conditions when
 minimizing a quadratic function with Hessian $G \in R^{n,n}$ subject
to the linear equality constraint $Cp = c, \ 
   C \in {R}^{m,n}, \ p,g \in R^n, \ 
   c,z \in R^m $
 $$
 \left [
   \begin{array}{ll}
    G  & \,\, C^T\\
       &\\
    C  & \, \, 0
    \end{array}
    \right ] 
    \left (
    \begin{array}{l}
    p \\ z
    \end{array}
    \right )    =
    \left (
    \begin{array}{l}
    g \\ c
    \end{array}
    \right )
\eqno{(15)}
$$

If $G$ is nonsingular, the coefficient matrix is nonsingular iff $CG^{-1}C^T$
is nonsingular. Usually $G$ is nonsingular, symmetric and positive definite,
but this assumption, required by several classical solvers, is not
necessary for  the ABS solvers. 

  To derive ABS methods using the structure of system (15), observe that 
(15) is equivalent to the two subsystems
$$
Gp + C^Tz = g    \eqno{(16)}
$$
$$
Cp = c.
\eqno{(17)}
$$
Consider the general solution of $Cp=c$  in the  ABS form, with
$q \in R^n$ arbitrary
$$
p = p_{m+1}+ H_{m+1}^Tq
    \eqno{(18)}
$$
The parameters used to construct $p_{m+1}$
and $H_{m+1}$ are arbitrary, hence (18) defines a class of algorithms.

Since the KT equations have a unique solution, there is a
 $q$  which makes $p$  the unique   
$n$-dimensional subvector defined by the first $n$ components of
the solution of (15).
By multiplying  $Gp + C^Tz = g$ on the left by $H_{m+1}$
  we obtain the equation
$$
H_{m+1}Gp = H_{m+1}g    \eqno{(19)}
$$
which does not contain $z$. Now there are two possibilities for  determining
$p$:
\begin{description}
\item[(A1)]
Consider the system formed by  (19) and (17). Such a system
is solvable but overdetermined.  Since  $\mbox{rank}(H_{m+1}) =n-m$,
$m$ equations are recognized as dependent and
are eliminated in step (B) of any ABS algorithm applied to this system,
which then computes the unique solution.
\item[(A2)]
In equation (19)  replace $p$ by  the  general
solution (18) to give
$$
H_{m+1}GH_{m+1}^Tq=H_{m+1}g-H_{m+1}Gp_{m+1}.
\eqno{(20)}
$$
The above system can be solved by any ABS method for a
particular solution $q$, $m$ equations being again removed at step {\bf (B)}
of  the  ABS
algorithm as linearly dependent. 
\end{description}
Once $p$ is determined, one can determine
$z$ in two ways, namely:
\begin{description}
\item[(B1)]
Solve by any ABS method the overdetermined compatible system
$$
C^Tz = g - Gp
\eqno{(21)}
$$
by removing at step (B) of the ABS algorithm the $n-m$ dependent equations.\\
\item[(B2)]
Let $P = (p_1,\ldots, p_m)$ be the matrix whose columns are the search vectors
generated on the system $Cp = c$. Now  $CP = L$, with $L$
nonsingular lower diagonal. Multiplying equation (21) on the left by
$P^T$ we obtain a
triangular system,  defining $z$ uniquely
$$
L^Tz = P^Tg - P^TGp.
\eqno{(22)}
$$
\end{description}
Extensive numerical testing has  evaluated the accuracy
of the above  ABS algorithms for KT equations for certain choices
of the ABS parameters (corresponding to the 
implicit LU algorithm with row pivoting
and the modified Huang algorithm). The methods have been tested against
 the method of Aasen and methods using the LU and the
QR factorization. The experiments have shown that some ABS methods are
the most accurate, in both residual and solution error; moreover some ABS
algorithms are cheaper in storage and in overhead, up to one order, 
especially for the case when $m$ is close to $n$. In particular two methods
based upon the implicit LU algorithm  not only have turned out to be more accurate, especially
in residual error, than the method of Aasen and the method using QR
factorization via Houselder matrices, but  are also cheaper in number
of operations (the method of Aasen has a lower storage for small $m$ but
a higher storage for large $m$).

In many interior point methods the main computational cost is to compute the
solution for a sequence of KT problems where only $G$, which is diagonal, changes.
In such a case the ABS methods, which initially work on the matrix $C$,
which is unchanged, have  an advantage, particularly when
$m$ is large, where the dominant cubic term decreases with $m$ and disappears for
$m=n$, so that the overhead is dominated by  second order terms. Again
numerical experiments show that some ABS methods are more accurate than
the classical ones. 
\vskip9mm
\section{A class of ABS methods for matrix equations}
\vskip5mm
It is common, in particular in optimization, to find  systems of
matrix equations
$$
A^i \bullet X=b_i, \quad i=1, \ldots, m
\eqno(23)
$$
where operation $\bullet$ is defined by $A \bullet B=$tr$(A^TB)$. We can
write systems (23), with obvious definition of $\circ$, as
$$
 {\cal A} \circ X =b
 \eqno{(24)}
$$
where ${\cal A}$ has the following form, with 
 $A^{i} \in {R}^{n,n}$ for $i=1, \ldots,m$
$$
{\cal A}=\left
	 [
	 \begin{array}{l}
	 A^{1}\\
	 \vdots\\
	 A^{m}
	 \end{array}
	 \right
	 ].       \eqno{(25)}
$$

\noindent Problem (24) is a linear system in the space of matrices and the
associated projection operators are
matrices whose elements are matrices in ${R}^{n,n}$.  This
observation led us to consider a linear space denoted by
$({R}^{n,n})^{n,n}$  and study
its linear algebra.
The isomorphisms between
${R}^{n,n}$ and ${R}^{n^2}$, between
$({R}^{n,n})^{n,n}$ and ${R}^{n^2,n^2}$ allow to
establish  ABS algorithms to solve (24), to
 generalize the
 Huang algorithm and implicit LU algorithm and to define other special
algorithms.
Quasi-Newton matrices satisfying
 linear relations
(for example, symmetry and sparsity) can be described in the considered matrix 
form
 and can be solved by
the proposed matrix ABS algorithm. 

The ABS method for finding a solution
$X \in {R}^{n,n}$, of  system
(24) is as follows, the symbol $^*$ indicating transposition in the
matrix space.
\vskip3mm
\noindent {\bf The matrix ABS algorithm}

\begin{description}
\item[Step 1 \mbox{}]Give $X^{1} \in {R^{n,n}}^{n,n},
{\cal H}^{1} \in ({R})^{n,n}$, set $k=1.$
\item[Step 2 \mbox{}]Compute $\tau_k=A^k \bullet X^k-b_k$ and $S^k=
{\cal H}^k \circ A^k$.
\item[Step 3 \mbox{}]If $S^k \neq 0$ go to Step 4; if $S^k=0$ and
$\tau = 0$ set $X^{k+1}=X^k,\ {\cal H}^{k+1}={\cal H}^k$ and go to Step 7
if $k< m$; otherwise stop.  If $S^k =0$ and $\tau \neq 0$ stop, the system
(24) is incompatible.
\item[Step 4 \mbox{}]Compute $P^k \in {R}^{n,n}$ by
$$
P^k=({\cal H}^k)^* \circ Z^k
\eqno{(26)}
$$
where $Z^k \in {R}^{n,n}$ is arbitrary save that
$
A^k \bullet P^k \neq 0
$
\item[Step 5 \mbox{}]Update the approximation of the solution by
$$
X^{k+1}=X^k-\alpha_k P^k, \ \ \alpha_k = \tau_k/A^k \bullet P^k
\eqno{(27)}
$$
If $k=m$ stop; $X^{m+1}$ solves the system.
\item[Step 6 \mbox{}]Update the matrix ${\cal H}^k$ by
$$
{\cal H}^{k+1}={\cal H}^k-{\cal H}^k \circ A^k \otimes [({\cal H}^k)^*
\bullet W^k]/W^k \bullet ({\cal H}^k \circ A^k)
\eqno{(28)}
$$
where $W^k \in {R}^{n,n}$ is arbitrary save that
$
W^k \bullet ({\cal H}^k \circ A^k) \neq 0
$
\item[Step 7 \mbox{}]Increment the index $k$ by one and goto Step 2.
\end{description}
 Properties of the  matrix ABS method  generalize properties
 of the scaled ABS class, albeit proofs are not always obvious, 
 see Spedicato, Xia and Zhang (1999). We just
 recall that
the general solution of (24) can be expressed as follows, with  
$W \in {R}^{n,n}$ arbitrary
$$
X =X^{m+1}+
({\cal H}^{m+1})^* \circ W
\eqno{(29)}
$$
There are natural generalizations of the Huang and the implicit LU algorithms
in the ABS class.
 Huang matrix 
algorithm allows to construct  solutions to
the quasi-Newton equation
$$
B' \delta = r
\eqno{(30)}
$$
where $\delta=x'-x$ and $r=F(x')-F(x)$ when solving the nonlinear system of
equations $F(x)=0$, $r=\nabla f(x')- \nabla f(x)$ when minimizing the
unconstrained function
$ f(x)$.
Equation (30) can be solved also under the
additional condition that some elements of the solution take prescribed values,
by which way we can introduce the conditions of sparsity, symmetry 
and positive definiteness.
\vskip9mm
\section{A class of ABS derived iterative methods}
\vskip5mm
\noindent When $n$ is large, both storage  and number of operations
may be too large for an implementation of  ABS methods using
explicitly the Abaffians. One can develop methods with
lower storage and 
operations  by working with  formulations of ABS methods
that use vectors and then using restart or truncation.
This approach leads  to loss of termination and
generates iterative methods. Here we  describe one such class, which
is derived from the following formulation of ABS methods in terms
of $k$ vectors at the $k-$th step, due originally to Bodon, see
Abaffy and Spedicato (1989), and where we take $z_k, w_k$ to be multiple of
each other.
\vskip3mm
\noindent {\bf The Bodon-ABS vector algorithm}

\vskip10pt

$\quad $ {\bf Let} $x_1\in R^n$ { be arbitrary}. {Let}
$H_1\in R^{n,n}$ and $V=(v_1,\ldots,v_n)\in R^{n,n}$
 be arbitrary

$\quad$ nonsingular.

$\quad$ {\bf For} $k=1$ to $n$

$\quad\qquad\tau_k=v^T_kAx_k-v^T_kb$

$\quad\qquad$ {\bf If} $k>1$ then

$\quad\qquad\qquad p^1_k=H^T_1z_k$

$\quad\qquad\qquad$ {\bf For } $j=1$ to $k-1$

$\quad\qquad\qquad\qquad
p^{j+1}_k=p^j_k-\left({v^T_jAp^j_k}/{v^T_jAp_j}\right)p_j
$

$\quad\qquad\qquad$ {\bf End}

$\qquad\quad\qquad p_k=p_k^k$

$\quad\qquad$ {\bf Else}

$\quad\qquad\qquad$ $p_k=H^T_1z_1$

$\quad\qquad$ {\bf Endif}

$\quad\qquad\alpha_k={\tau_k}/{(v^T_kAp_k)}$

$\quad\qquad
x_{k+1}=x_k-\alpha_kp_k
$

$\quad${\bf End}

\vskip20pt

\noindent The above algorithm leads to iterative methods either via restart or via
truncation. Here we consider only truncation.
If $m$ is the number of available storage vectors,
then we keep only information from the last $m$
iterations, i.e.  we replace all iterations from $j=1$ to $k-1$ by
iterations from $j=k-m$ to $k-1$.
Strategies where the kept vectors are not necessarily the last $m$ vectors
may of course be considered.
The matrix
$H_1$ should  require  low storage, hence we take
 $H_1=I$. Parameters $v_k$ should be linearly
 independent and parameters $z_k$  feasible (no division by zero).
\vskip3mm

 \noindent{\bf Algorithm ABS(m): the truncated Bodon-ABS vector algorithm}

\vskip10pt

$\quad $ {\bf Let} $x_1\in R^n$ { be arbitrary}. 

$\quad$ Give the integer $m, 1\le m\le n$.

$\quad$ {\bf Do} $k=1,2,\cdots$, until convergence

$\quad\qquad\tau_k=v^T_kAx_k-v^T_kb$

$\quad\qquad$ {\bf If} $k>1$ then

$\quad\qquad\qquad t=\max(1,k-m)$

$\quad\qquad\qquad p^t_k=z_k$

$\quad\qquad\qquad$ {\bf For } $j=t$ to $k-1$

$\quad\qquad\qquad\qquad
p^{j+1}_k=p^j_k-\left({v^T_jAp^j_k}/{v^T_jAp_j}\right)p_j
$

$\quad\qquad\qquad$ {\bf End}

$\qquad\quad\qquad p_k=p_k^k$

$\quad\qquad$ {\bf Else}

$\quad\qquad\qquad$ $p_k=z_1, \ (z_1\in R^n)$

$\quad\qquad$ {\bf Endif}

$\quad\qquad\alpha_k={\tau_k}/{(v^T_kAp_k)}$

$\quad\qquad
x_{k+1}=x_k-\alpha_kp_k
$

$\quad${\bf End do}

\vskip20pt

If $m=1$ one can show that  the truncated
implicit LU, 
 Huang and implicit QR algorithm generate
 respectively  the Gauss-Seidel, the Kaczmarz
and the De la Garza methods. For a comparison of storage and operations
requirements with other well-know iterative methods see 
Spedicato and Li (1999).

One can prove that  algorithm ABS($m$) is well-defined
if, letting $V_k=(v_{t},\cdots,v_{k})$ and
$Z_k=(z_{t},\cdots,z_{k})$, where $t=\max(1,k-m)$, then $V^T_kAZ_k$
is strongly nonsingular
for all $k$.

Without loss of generality, 
we can define the scaling
parameters by
$$
v_k=A^{-T}Yp_k \eqno(31)
$$
where $Y$ is a symmetric, positive definite matrix.
The following choices for $Y$ define, in the original scaled 
ABS class, the three   subclasses considered in section 1
and require a low storage for large sparse systems: 
$Y=I$ (the optimally scaled subclass); 
$Y=A^TA$ (the orthogonally scaled subclass);
$Y=A^T>0$, for $A=A^T>0$, (the conjugate direction subclass).

 Algorithm ABS$(m)$
with the parameter choice
(31) has   variational
properties, related to those of the original ABS class, namely
 that $x_{k+1}$ minimizes the error
 $Y$-weighted  Euclidean norm over a linear variety spanned by the last 
 $m$ search vectors.

One can also  show that the truncated
generalized conjugate
direction algorithm of Dennis and Turner (1987) can be
obtained by
special choices of the parameters in Algorithm ABS($m$). From this
equivalence  it follows  that 
ORTHOMIN$(m))$,
ORTHODIR$(m))$ and
 GMRES($m)$  are special cases of ABS($m$).

The convergence of  ABS($m$) at a linear rate  can be proved
 by requiring
 that the angle between 
 $p_k$ and the gradient $\overline r_{k}$ 
of the $Y-$weighted error norm be uniformly bounded away from $90^0$.
\vskip3mm
\begin{Theo}. \   
{\it Suppose that
 there exists $\gamma>0$ such that for the search vectors generated
by Algorithm ABS($m$) with choice (31) one has
$$\mid p^T_k\overline
r_{k}\mid\ge \gamma\parallel p_k\parallel_2\parallel
\overline r_{k}\parallel_2\eqno(32)
$$
for all $k$. Then the sequence $x_k$
converges to $x^\ast$ and satisfies
$$
\parallel
x_{k+1}-x^\ast\parallel_Y\le
(1-\frac{\gamma^2}{Cond(Y)})^{1/2}
\parallel x_{k}-x^\ast\parallel_Y,\eqno(33)
$$
where $Cond(Y)$ is  defined by the ratio of the largest to the  smallest
eigenvalues of $Y$.}
\end{Theo}

\vskip9mm
\section{ABSPACK and its numerical performance}
\vskip5mm
The ABSPACK  project aims at producing a mathematical package for solving
linear and nonlinear systems and optimization problems using the ABS
algorithms. The project will take several years for completion, in view
of the substantial work needed to test  the  alternative ways
ABS methods can be implemented (via different linear algebra formulations
of the process,  different possibilities of reprojections,  different
possible block formulations etc.) and of the necessity of comparing  the
produced software with the established packages in the market (e.g. 
MATLAB, LINPACK, LAPACK,  UFO ...). It is expected that the software
will be documented in a forthcoming monograph and willl be made available
to general users.

At the present state of the work FORTRAN 77 implementations have been made
of several versions of the following ABS algorithms for  solving  linear
systems:
\begin{enumerate}
\item  The Huang  and the modified Huang algorithms in two different
       linear algebra versions of the process, for solving determined,
	underdetermined  and overdetermined systems,
       for  a solution of least Euclidean norm
\item  The implicit LU and implicit LX algorithm for determined, underdetermined
       and overdetermined linear systems, for a solution of basic type
\item  The implicit QR algorithm for determined, underdetermined and
       overdetermined linear systems, for a solution of basic type
\item  The above algorithms for some structured problems (including
       banded and block angular matrices)
\item  An ABS algorithm proposed by Adib and  Mahdavi-Amiri (1999) equivalent
       to a block-two ABS method
\item  ABS versions of the GMRES method, some requiring  less storage
\item  Several ABS methods  for KT equations.
\end{enumerate}
For  a full presentation of the above methods  and their comparison with
NAG, LAPACK, LINPACK and UFO codes see Bodon, Luksan and Spedicato (2000),
Bodon and Spedicato (2000a,b). Some results are presented in the Appendix.
There the  columns refer respectively to: 
the problem, the dimension, the  algorithm, the relative solution error
(in Euclidean norm), the relative residual error in Euclidean norm
(i.e. ratio of residual
norm over norm of right hand side), the computed rank and the time in
seconds. Computations have been performed in double precision on  a
Digital Alpha workstation with machine zero about $10^{-17}$. All test
problems have been generated with integer entries or powers of two such
that all entries are exactly represented in the machine and 
the right hand side can be computed exactly, so that the given solution
is an exact solution of the problem as it is represented in the machine.
Comparison is given with some LAPACK and LINPACK codes, including those based
upon singular value decomposition ({\it svd}) and rank revealing QR
factorization  ({\it gqr}).

Analysis of all obtained results indicates:
\begin{enumerate}
\item  Modified Huang is generally the most accurate ABS algorithm and
compares in accuracy with the best LAPACK solvers based upon singular value
factorization and rank revealing QR factorization; also the estimated ranks
are usually the same.
\item  On problems whose numerical estimated rank is much less than the
dimension, one of the versions  of modified Huang is much faster than the
LAPACK codes using SVD or rank revealing QR factorization, even more than
a factor 100. This is due to the fact 
that once an equation is recognized as
dependent it does not contribute to the general overhead in ABS algorithms.
\item Modified Huang is generally faster and more accurate than other
ABS methods and classical methods on KT equations.
\end{enumerate}
It should be noted that the performance of the considered ABS algorithms in
term of times  could be improved by developing block versions,  as it is
the case for the LAPACK codes, a work presently in progress.
\newpage
\section{Final remarks}
\vskip5mm
Additional work done in the framework of the project, not described here
in detail, has involved the following topics.
\begin{itemize}
\item Improvement of the performance of Newton method via a special
truncation approach, see Deng and Wang (1998).
\item Critical review  of variable metric methods for unconstrained
optimization with discussion of the ABS applications in this field, see
Luksan and Spedicato (1998).
\item Analysis of the relations between the ABS methods and the classical
method of  averaging functional corrections, see Gredzhuk and
Petrina (1998).
\item Development of indefinitely preconditioned truncated Newton methods for
large sparse equality constrained nonlinear programming problems, see
Luksan  and Vlcek (1998).
\item Computation and update of inertias of KKT matrices for use in  
quadratic programming, see Zhang (1999).
\item Further applications of the implicit LX algorithm to the simplex
method for the LP problem, see Spedicato and Xia (1999).
\end {itemize}
The field of ABS methods is now mature from a theoretical point of view,
albeit there are exciting possibilities for applications to new fields,
e.g. the eigenvalue problems. We expect that the completion of the
project ABSPACK will provide a useful new instrument for users of
mathematical software.

\newpage
\begin{verbatim}   
    APPENDIX


    RESULTS ON DETERMINED  LINEAR SYSTEMS
 
    Condition  number: 0.21D+20
    IDF2    2000     huang2       0.10D+01  0.69D-11    2000   262.00
    IDF2    2000     mod.huang2   0.14D+01  0.96D-12       4     7.00
    IDF2    2000     lu lapack    0.67D+04  0.18D-11    2000    53.00
    IDF2    2000     qr lapack    0.34D+04  0.92D-12    2000   137.00
    IDF2    2000     gqr lapack   0.10D+01  0.20D-14       3   226.00
    IDF2    2000     lu linpack   0.67D+04  0.18D-11    2000   136.00

    Condition number: 0.10D+61
    IR50    1000     huang2       0.46D+00  0.33D-09    1000    36.00
    IR50    1000     mod.huang2   0.46D+00  0.27D-14     772    61.00
    IR50    1000     lu lapack    0.12D+04  0.12D+04     972     7.00
    IR50    1000     qr lapack    0.63D+02  0.17D-12    1000    17.00
    IR50    1000     gqr lapack   0.46D+00  0.42D-14     772    29.00
    IR50    1000     lu linpack   --- break-down ---                 


    RESULTS ON OVERDETERMINED SYSTEMS
 
    Condition  number: 0.16D+21
    IDF3    1050  950    huang7       0.32D+04  0.52D-13    950    31.00
    IDF3    1050  950    mod.huang7   0.14D+04  0.20D-09      2     0.00
    IDF3    1050  950    qr lapack    0.37D+13  0.83D-02    950    17.00
    IDF3    1050  950    svd lapack   0.10D+01  0.24D-14      2   145.00
    IDF3    1050  950    gqr lapack   0.10D+01  0.22D-14      2    27.00


    Condition number:  0.63D+19
    IDF3    2000  400    huang7       0.38D+04  0.35D-12    400     9.00
    IDF3    2000  400    mod.huang7   0.44D+03  0.67D-12      2     0.00
    IDF3    2000  400    impl.qr5     0.44D+03  0.62D-16      2     0.00
    IDF3    2000  400    expl.qr      0.10D+01  0.62D-03      2     0.00
    IDF3    2000  400    qr lapack    0.45D+12  0.24D-02    400     8.00
    IDF3    2000  400    svd lapack   0.10D+01  0.65D-15      2    17.00
    IDF3    2000  400    gqr lapack   0.10D+01  0.19D-14      2    12.00
    
    RESULTS ON UNDERDETERMINED LINEAR SYSTEMS

    Condition number:  0.29D+18
    IDF2     400 2000    huang2       0.12D-10  0.10D-12    400    12.00
    IDF2     400 2000    mod.huang2   0.36D-08  0.61D-10      3     1.00
    IDF2     400 2000    qr lapack    0.29D+03  0.37D-14    400     9.00
    IDF2     400 2000    svd lapack   0.43D-13  0.22D-14      3    68.00
    IDF2     400 2000    gqr lapack   0.18D-13  0.24D-14      3    12.00
  
    Condition number:  0.24D+19
    IDF3     950 1050    huang2       0.00D+00  0.00D+00    950    33.00
    IDF3     950 1050    mod.huang2   0.00D+00  0.00D+00      2     1.00
    IDF3     950 1050    qr lapack    0.24D+03  0.56D-14    950    17.00
    IDF3     950 1050    svd lapack   0.17D-14  0.92D-16      2   178.00
    IDF3     950 1050    gqr lapack   0.21D-14  0.55D-15      2    26.00
\end{verbatim}
\newpage
\begin{verbatim}

    RESULTS ON KT SYSTEMS

    Condition number:  0.26D+21 
    IDF2    1000  900    mod.huang    0.55D+01  0.23D-14     16    24.00
    IDF2    1000  900    impl.lu8     0.44D+13  0.21D-03   1900    18.00
    IDF2    1000  900    impl.lu9     0.12D+15  0.80D-02   1900    21.00
    IDF2    1000  900    lu lapack    0.25D+03  0.31D-13   1900    62.00
    IDF2    1000  900    range space  0.16D+05  0.14D-11   1900    87.00
    IDF2    1000  900    null space   0.89D+03  0.15D-12   1900    93.00
 
 
    Condition number:  0.70D+20 
    IDF2    1200  600    mod.huang    0.62D+01  0.20D-14     17    36.00
    IDF2    1200  600    impl.lu8     0.22D+07  0.10D-08   1800    44.00
    IDF2    1200  600    impl.lu9     0.21D+06  0.56D-09   1800    33.00
    IDF2    1200  600    lu lapack    0.10D+03  0.79D-14   1800    47.00
    IDF2    1200  600    range space  0.11D+05  0.15D-11   1800    63.00
    IDF2    1200  600    null space   0.38D+04  0.13D-12   1800   105.00
\end{verbatim}
\end{document}